\documentclass[onecolumn]{revtex4}

\usepackage{amssymb,amsmath}

\usepackage{natbib}

\newtheorem{definition}{Definition}

\begin{document}

\title{The regularized product of the Fibonacci numbers}
\author{Adrian R. Kitson}
\email{a.r.kitson@massey.ac.nz}
\affiliation{Institute of Fundamental Sciences, Massey University, Private Bag 11 222, Palmerston North, New Zealand}

\maketitle

Using zeta-function regularization, it is possible to calculate, for example, the regularized product over all the positive integers,
\begin{displaymath}
\prod_{n=1}^{\infty}n=\sqrt{2\pi},
\end{displaymath}
and over all the primes~\cite{primes},
\begin{displaymath}
\prod_{p}p=\left(2\pi\right)^{2}.
\end{displaymath}
What is the regularized product over all the Fibonacci numbers? In other words, what is:
\begin{displaymath}
\Delta=\prod_{n=1}^{\infty}f_{n}^{\mbox{}},
\end{displaymath}
where \(f_{n}^{\mbox{}}\) is the \(n\)th Fibonacci number? First we define what we mean by `regularized product.' 
\begin{definition}
The \emph{regularized product} of \(0<a_{1}^{\mbox{}}\le a_{2}^{\mbox{}}\le a_{3}^{\mbox{}}\le\dots\) is
\begin{displaymath}
\prod_{n=1}^{\infty} a_{n}^{\mbox{}}=\exp\left(-\zeta_{A}^{\prime}\left(0\right)\right),
\end{displaymath}
where \(\zeta_{A}\) is the regularized zeta function
\begin{displaymath}
\zeta_{A}^{\mbox{}}\left(s\right)=\sum_{n=1}^{\infty}a_{n}^{-s}.
\end{displaymath}
\end{definition}
The \emph{Fibonacci numbers} are \(f_{1}^{\mbox{}}=1\), \(f_{2}^{\mbox{}}=1\) and then
\begin{displaymath}
f_{n}^{\mbox{}}=f_{n-1}^{\mbox{}}+f_{n-2}^{\mbox{}},
\end{displaymath}
for \(n\ge 2\). The \(n\)th Fibonacci number can be expressed as
\begin{displaymath}
f_{n}^{\mbox{}}=\frac{\phi^{n}-\left(-\phi\right)^{-n}}{\sqrt{5}},
\end{displaymath}
where \(\phi=\left(1+\sqrt{5}\,\right)\!/2\) is the \emph{Golden Mean}. The Fibonacci zeta function is
\begin{displaymath}
\zeta_{F}^{\mbox{}}\left(s\right)=\sum_{n=1}^{\infty}f_{n}^{-s},
\end{displaymath}
which, by the ratio test, is absolutely convergent for \(Re\left(s\right)>0\). Its derivative with respect to \(s\) is
\begin{displaymath}
\zeta_{F}^{\prime}\left(s\right)=-\sum_{n=1}^{\infty}\ln\left(f_{n}^{\mbox{}}\right)f_{n}^{-s}.
\end{displaymath}
From the definition, the regularized product of the Fibonacci numbers is
\begin{equation}
\label{eq:delta}
\Delta=\exp\left(-\zeta_{F}^{\prime}\left(0\right)\right).
\end{equation}
Unfortunately, the derivative of the Fibonacci zeta function diverges as \(s\rightarrow 0\) and therefore must be regularized. Adding and subtracting the leading divergence gives
\begin{align}
\zeta_{F}^{\prime}\left(s\right)&=-\sum_{n=1}^{\infty}\ln\left(\frac{\phi^{n}}{\sqrt{5}}\right)\left(\frac{\phi^{n}}{\sqrt{5}}\right)^{-s}-\sum_{n=1}^{\infty}\left[\ln\left(f_{n}^{\mbox{}}\right)f_{n}^{-s}-\ln\left(\frac{\phi^{n}}{\sqrt{5}}\right)\left(\frac{\phi^{n}}{\sqrt{5}}\right)^{-s}\right],\nonumber\\
&=\frac{5^{s/2}\left(\left(\phi^{s}-1\right)\ln\left(5\right)-2\phi^{s}\ln\left(\phi\right)\right)}{2\left(\phi^{s}-1\right)^{2}}-\sum_{n=1}^{\infty}\left[\ln\left(f_{n}^{\mbox{}}\right)f_{n}^{-s}-\ln\left(\frac{\phi^{n}}{\sqrt{5}}\right)\left(\frac{\phi^{n}}{\sqrt{5}}\right)^{-s}\right].\nonumber
\end{align}
The first term is a meromorphic function and the second term converges as \(s\rightarrow 0\). The Laurent expansion about \(s=0\) is
\begin{displaymath}
\zeta_{F}^{\prime}\left(s\right)=-\frac{1}{\ln\left(\phi\right)s^{2}}+\frac{1}{24}\left(2\ln\left(\phi\right)+\frac{3\ln^{2}\left(5\right)}{\ln\left(\phi\right)}-6\ln\left(5\right)\right)-\sum_{n=1}^{\infty}\left[\ln\left(f_{n}^{\mbox{}}\right)-\ln\left(\frac{\phi^{n}}{\sqrt{5}}\right)\right]+{\mathcal O}\left(s\right).\nonumber
\end{displaymath}
The sum simplifies,
\begin{align}
\sum_{n=1}^{\infty}\left[\ln\left(f_{n}^{\mbox{}}\right)-\ln\left(\frac{\phi^{n}}{\sqrt{5}}\right)\right]&=\sum_{n=1}^{\infty}\ln\left(1-\left(-\phi\right)^{-2n}\right),\nonumber\\
&=\ln\left(\prod_{n=1}^{\infty}\left(1-\left(-\phi\right)^{-2n}\right)\right),\nonumber\\
&=\ln\left(c\right),\nonumber
\end{align}
where \(c\) is the \emph{Fibonacci factorial constant}~\cite{finch,sloane}. Thus,
\begin{displaymath}
\zeta_{F}^{\prime}\left(s\right)=-\frac{1}{\ln\left(\phi\right)s^{2}}+\frac{1}{24}\left(2\ln\left(\phi\right)+\frac{3\ln^{2}\left(5\right)}{\ln\left(\phi\right)}-6\ln\left(5\right)\right)-\ln\left(c\right)+{\mathcal O}\left(s\right).
\end{displaymath}
The regularized value is the limit as \(s\rightarrow 0\) once the principal part has been removed, that is,
\begin{displaymath}
\zeta_{F}^{\prime}\left(0\right)=\frac{1}{24}\left(2\ln\left(\phi\right)+\frac{3\ln^{2}\left(5\right)}{\ln\left(\phi\right)}-6\ln\left(5\right)\right)-\ln\left(c\right).
\end{displaymath}
Using equation~(\ref{eq:delta}),
\begin{displaymath}
\Delta=\frac{5^{1/4}\exp\left(-\ln^{2}\left(5\right)/8\ln\left(\phi\right)\right)c}{\phi^{1/12}}.
\end{displaymath}
The Fibonacci factorial constant can be expressed in terms of the derivative of the Jacobi theta function of the first kind. Consequently,
\begin{displaymath}
\Delta=\left(-1\right)^{1/24}5^{1/4}\exp\left(-\ln^{2}\left(5\right)/8\ln\left(\phi\right)\right)\left(\frac{1}{2}\vartheta_{1}^{\prime}\left(0,-\frac{i}{\phi}\right)\right)^{1/3}.
\end{displaymath}
Numerically, \(\Delta=0.8992126807\ldots\).

\begin{acknowledgments}
Marie-Claire Kitson is acknowledged for useful conversations.
\end{acknowledgments}

\end{document}